\documentclass[12pt]{amsart}
\usepackage{amsmath,amsfonts,amssymb,amscd}
\usepackage{latexsym,epsfig}
\usepackage[latin1]{inputenc}
\usepackage[T1]{fontenc}
\usepackage{enumerate}
\usepackage{amsbsy}
\usepackage[matrix,arrow]{xy}

\newcommand{\N}{\ensuremath{\mathbb N}}
\newcommand{\Z}{\ensuremath{\mathbb Z}}
\newtheorem{thm}{Theorem}
\newtheorem{lem}[thm]{Lemma}
\newtheorem{prop}[thm]{Proposition}
\newtheorem{cor}[thm]{Corollary}
\newtheorem{rem}[thm]{Remark}

\begin{document}

\title{Combinatorial properties of virtual braids}

\author[Bardakov]{Valerij G. Bardakov}
\address{Sobolev Institute of Mathematics, Novosibirsk 630090, Russia}
\email{bardakov@math.nsc.ru}

\author[Bellingeri]{Paolo Bellingeri}
\address{Univ. Nantes, Laboratoire J.~Leray, 44322 Nantes, France}
\email{paolo.bellingeri@univ-nantes.fr}





\begin{abstract}
We study combinatorial properties of virtual braid groups and we
describe relations with finite type invariant theory for virtual
knots and Yang Baxter equation.
\end{abstract}

\maketitle

\section{introduction}

Virtual knots were introduced by L.~Kauffman~\cite{Ka} as the geometric counterpart of
Gauss diagrams.
A virtual knot diagram is a generic (oriented) immersion of the circle into the plane,
with the usual positive and negative
crossings plus a new kind of crossings called virtual. The Gauss diagram of a virtual
knot is constructed in the same way
as for a classical knot excepted that virtual crossings are disregarded.
Gauss diagrams turn to be in bijection
with virtual knots diagrams up to isotopy and a finite number of  ``virtual''
 moves around crossings, which generalise
usual Reidemeister moves.
Using virtual Reidemester moves we can introduce  a notion of ``virtual'' braids
(see for instance
\cite{B,K,Ver01}). Virtual braids on $n$ strands form a group, usually denoted by $VB_n$.
The relations between virtual braids and virtual knots (and links) are completely
determined by a generalisation
of Alexander and  Markov Theorems~\cite{K}.

The paper is mainly devoted to the study of the kernels of two different projections
of $VB_n$ in $S_n$,
the normal closure of the braid group $B_n$, that we will denote by $H_n$, and
the so-called virtual pure braid group
$VP_n$, which is related to the quantum Yang Baxter
equation (see Section~\ref{vp}).

The paper is organised as follows:
in Section~\ref{intro} we recall some definitions and classical
results on combinatorial group theory and in Section~\ref{braids}
we provide a short survey on lower central series for
(generalised) braids and their  relations with finite type
invariant theory. In following Sections we determine the lower
central series of the virtual braid group $VB_n$ and of its
subgroups $H_n$ and $VP_n$. Finally, in Section~7 we determine
the intersection of $H_n$ and $VP_n$ and in  Section~\ref{gpv} we
provide a connection between virtual pure braids and  the finite
type invariant theory for virtual knots defined by Goussarov,
Polyak and Viro in~\cite{GPV}.

\section{Definitions} \label{intro}

Given  a group $G$, we define the \emph{lower central series} of $G$
as the filtration  $\Gamma_1(G)=G \supseteq \Gamma_2(G) \supseteq \ldots$, where
$\Gamma_i(G)=[\Gamma_{i-1}(G), G]$.
 The \emph{rational lower central series}
of $G$ is the   filtration $D_1(G) \supseteq  D_2(G)
\supseteq \ldots$  obtained setting $D_1(G)=G$, and for $i\geq 2$,
defining $D_i(G)=\{ \, x \in G \, | \, x^n \in \Gamma_{i}(G) $ for some
$n\in \N^* \, \}$. This filtration was first considered by Stallings \cite{St} and
we denote it by the name proposed in~\cite{GL}.

Let us recall that for any
group-theoretic property $\mathcal{P}$, a group $G$ is said to be
\emph{residually $\mathcal{P}$} if for any (non-trivial) element $x
\in G$, there exists a group $H$ with the property  $\mathcal{P}$ and
a surjective homomorphism $\phi: G \to H$ such that  $\phi(x) \not=1$. It is well
known that a group $G$ is residually nilpotent if and only if
$\bigcap_{i \ge 1}\Gamma_i(G)=\{ 1\}$. On the other hand, a group $G$
is residually torsion-free nilpotent if and only if
$\bigcap_{i \ge 1}D_i(G)=\{ 1\}$.

The fact that a group is residually torsion-free nilpotent has several
important consequences, notably that the group is bi-orderable~\cite{MR}. We
recall that a group $G$ is said to be
\emph{bi-orderable} if there exists a strict total ordering $<$ on its
elements which is invariant under left and right multiplication, in
other words, $g<h$ implies that $gk<hk$ and  $kg<kh$ for all $g,h,k\in
G$.

Finally we recall that the  \emph{augmentation ideal} of a group $G$ is the
two-sided ideal $I(G)$ of the group ring $\Z[G]$ generated by  the set $\{g-1
\, | \, g\in G \} $. We denote by $I^d(G)$ the $d$th power of  $I(G)$.

The residually torsion free nilpotence of a group $G$ implies that
$\bigcap_{d \ge 1}I^d(G)=\{ 1\}$ (\cite{pas}, Theorem 2.15, chapter VI).

\section{Lower central series for generalized braid groups} \label{braids}

\noindent \textbf{Artin-Tits groups and surface braid groups.}
Let us start by recalling some standard results on combinatorial
properties of braid groups.  It is well known
(see~\cite{GoL} for instance) that  the Artin braid group $B_n$ is not residually
nilpotent for  $n\ge 3$ and that  the abelianization of $B_n$ is isomorphic to
$\Z$ and
$\Gamma_2(B_n)=\Gamma_3(B_n)$. We recall that classical braid groups are also called Artin-Tits
groups of type $\mathcal{A}$.  More precisely, let  $(W,S)$ be 
a Coxeter system and let us denote 
the order of the element $s t$ in $W$ by $m_{s, \, t}$ (for $s, t \in S$).
Let $B_W$ be the group defined by the following group presentation:
$$
B_W=\langle S \, | \underbrace{s t\cdots}_{m_{s, \, t}}= \underbrace{t s\cdots}_{m_{s, \, t}}\; \mbox{ for any $s\not=t \in S$
with $m_{s, \, t} < +\infty$} \, \rangle\;.
$$
The group $B_W$ is the Artin-Tits group associated to $W$. The group
$B_W$ is said to be of spherical type if  $W$ is finite. The kernel of
the canonical projection of $B_W$ onto $W$ is called the pure Artin-Tits group associated to $W$.
As explained in \cite{BGG},  it is easy to show
that the lower central series of
almost all Artin-Tits groups of spherical type also
stabilise at the second term, the only exception being Artin-Tits groups associated to the dihedral group
$I_{2m}$, for $m>1$.

When one considers surface braid groups (see \cite{BGG} for a definition)  new features appear. Let
$B_n(\Sigma)$ be the braid group on $n$ strands on the surface $\Sigma$.
In~\cite{BGG} it is proved that, when $\Sigma$ is an oriented surface of positive
genus and $n\ge 3$,  the  lower central series of $B_n(\Sigma)$ stabilises at
the third term.
Moreover, the quotient groups associated to the lower central series form a complete (abelian)
invariant for  braid groups of closed  surfaces.

\vspace*{5pt}

\noindent \textbf{Pure braids and finite type invariants.} Let $A, C$ be two groups.
If   $C$ acts on $A$ and the induced action on the abelianization
of $A$ is trivial, we say that $A \rtimes C$ is an
\emph{almost-direct product} of $A$ and $C$.

\begin{prop} {\bf (\cite{FR})}
Let $A, C$ be two groups. If   $C$ acts on $A$ and the induced action on the abelianization
of $A$ is trivial, then
$$
I(A \rtimes C)^d = \sum_{k=0}^{d} I(A)^k \otimes I(C)^{d-k} \quad \mbox{for all} \;  \;d \ge 0 \, .
$$
\end{prop}

\begin{prop}\label{corrn}
If $A \rtimes C$ is an \emph{almost-direct product} then
$\Gamma_m(A \rtimes C)=\Gamma_m(A) \rtimes \Gamma_m(C)$ and
$D_m(A \rtimes C)=D_m(A) \rtimes D_m(C).$
\end{prop}
\begin{proof}
The first statement is proved in \cite{FR2}.
Remark that the exact  sequence $1\to A \to A \rtimes C \to C \to 1$
 induces the following exact (splitting) sequence:
 $$
 1 \to D_m(A \rtimes C) \cap A \to D_m(A \rtimes C) \to  D_m(C) \to 1\,.
 $$
 Therefore the second statement is equivalent to prove that  $D_m(A \rtimes C) \cap A=
D_m(A)$ which is a straigthforward consequence
of the first statement and of the definition of rational lower central series.
\end{proof}

The structure of almost-direct product turns out to be a powerful tool in the determination
of algebras
related to lower central series (see for instance \cite{CCP}) and more generally in the
study of finite type invariants.

The pure braid group $P_n$ is an almost-direct product
of free groups~\cite{FR} and  this fact  has  been used in~\cite{pap}
in order to construct an universal finite type invariant
for braids with integers coefficients.
In~\cite{GP} Gonz\'alez-Meneses and Paris proved that
the  normal closure of the classical pure braid group $P_n$ in the pure braid group on $n$
strands of a surface $\Sigma$
is an almost-direct product of (infinitely generated) free groups. Adapting the approach of Papadima,
they constructed a universal finite type invariant for surface braids, but not multiplicative~\cite{BF}.

\begin{rem}
Free groups are residually torsion-free nilpotent~\cite{F,M}.
 It follows from Proposition~\ref{corrn} that pure braid groups are
residually torsion-free nilpotent~(see also \cite{FR2}). More generally,
using the faithfulness of the Krammer-Digne
representation, Marin has shown recently that the pure Artin-Tits
groups of spherical type are  residually torsion-free nilpotent~\cite{Ma}.
\end{rem}

\section{Lower central series of virtual braids}\label{virtual}

\begin{thm}{\bf(\cite{Ver01})} \label{presvb}
The group $VB_n$ admits the following group presentation:

\noindent {\bf $\bullet$ Generators:}  $\sigma_i$, $\rho_i$, $i =
1, 2, \ldots, n-1.$

\noindent {\bf $\bullet$ Relations:}
\begin{equation}
\sigma_i\sigma_{i+1}\sigma_i=\sigma_{i+1}
\sigma_i\sigma_{i+1},~~~i=1,2,\ldots,n-2 \,;  \label{eq2}
\end{equation}
\begin{equation}
\sigma_i\sigma_j=\sigma_j\sigma_i,~~~|i-j|\,   \geq 2 ;\label{eq1}
\end{equation}
\begin{equation}
\rho_{i}  \, \rho_{i+1} \, \rho_{i} = \rho_{i+1} \, \rho_{i} \, \rho_{i+1},~~~i=1,2,\ldots,n-2 \,;
\label{eq17}
\end{equation}
\begin{equation}
\rho_{i} \, \rho_{j} = \rho_{j} \, \rho_{i},~~~|i-j| \geq 2 \,;
\label{eq18}
\end{equation}
\begin{equation}
\rho_{i}^2 = 1,~~~i=1, 2, \ldots, n-1 \,;
\label{eq19}
\end{equation}
\begin{equation}
\sigma_{i} \, \rho_{j} = \rho_{j} \, \sigma_{i},~~~|i-j| \geq 2 \,;
\label{eq20}
\end{equation}
\begin{equation}
\rho_{i} \, \rho_{i+1}  \, \sigma_{i} = \sigma_{i+1} \, \rho_{i}  \, \rho_{i+1},
~~~i=1, 2, \ldots, n-2 \,.
\label{eq21}
\end{equation}
\end{thm}

\begin{rem}
Note that the last relation is equivalent to the following relation:
$$
\rho_{i+1} \, \rho_{i}  \, \sigma_{i+1} = \sigma_{i} \, \rho_{i+1}  \, \rho_{i}.
$$
On the other hand,  relations
$$
\rho_{i} \, \sigma_{i+1}  \, \sigma_{i} = \sigma_{i+1} \, \sigma_{i}
\, \rho_{i+1},~~\rho_{i+1} \, \sigma_{i}  \, \sigma_{i+1} = \sigma_{i}
\, \sigma_{i+1} \, \rho_i
$$
are not fulfilled in $VB_n$ (see for instance \cite{GPV}).
\end{rem}

\begin{prop}\label{embeds}
Let $\sigma_1, \dots, \sigma_{n-1}$ and respectively $s_1, \dots, s_{n-1}$
be the usual generators of the Artin braid group $B_n$ and of the symmetric group $S_n$.
The morphism $\iota: B_n \to VB_n$ defined by $\iota(\sigma_i)=\sigma_i$
and the morphism $\vartheta: S_n \to VB_n$ defined by $\vartheta(s_i)=\rho_i$
are well defined and injective.
\end{prop}
\begin{proof}
An easy argument for the injectivity of $\iota$ is given in \cite{K}.
Another different proof is given in \cite{R}.
Now, let $\mu: VB_n \to S_n$ be the morphism
defined as follows (this morphism will be considered in Section~\ref{defk}):
$$
\mu(\sigma_i)=1, \; \mu(\rho_i)=s_i, \;i=1,2,\dots, n-1\, ,
$$
The set-section $s: S_n \to VB_n$ defined by $\mu(s_i)=\rho_i$
is a well defined morphism and thus the subgroup generated by
$\rho_1, \dots, \rho_{n-1}$ is isomorphic to $S_n$ and
$\vartheta$ is injective.
\end{proof}

\noindent {\bf Notation.} In the following we use the notations $[a, b]= a^{-1} b^{-1} a b$
and $a^b=b^{-1} a b$.

\begin{prop} \label{braidlcs}
Let $VB_n$ be the virtual braid group on $n$ strands. The following properties hold:
\begin{enumerate}[a)]
\item The group $VB_2$ is isomorphic to $\Z*\Z_2$  which is
residually nilpotent.
\item The group $\Gamma_1(VB_n)/\Gamma_2(VB_n)$ is isomorphic to  $\Z\oplus\Z_2$ for $n \ge 2$.
\item The group $\Gamma_2(VB_3)/\Gamma_3(VB_3)$ is isomorphic to $\Z_2$. Otherwise, if  $n\ge 4$ then
 $\Gamma_2(VB_n)=\Gamma_3(VB_n)$.
\item If $n\ge 3$ the group  $VB_n$ is not residually nilpotent.
\end{enumerate}
\end{prop}
\begin{proof}
\begin{enumerate}[a)]
\item The group $\Z*\Z_2$ can be realised as a subgroup of $\Z_2*\Z_2*\Z_2$
  which is residually nilpotent (see~\cite{G} and~\cite{BGG}).
It  also follows from a result of Malcev~\cite{M}.
\item The statement can be easily verified
considering $VB_n$  provided with the group presentation given in  Theorem~\ref{presvb}.
\item
Consider the quotient $G = VB_3/\Gamma_3(VB_3)$. This group is
generated by $\overline{\sigma}_1 = \sigma_1 \Gamma_3(VB_3),$
$\overline{\sigma}_2 = \sigma_2 \Gamma_3(VB_3)$,
$\overline{\rho}_1 = \rho_1 \Gamma_3(VB_3)$ and $\overline{\rho}_2 = \rho_2 \Gamma_3(VB_3)$.
Since,
$$
\sigma_2 = [\sigma_1, \sigma_2] \sigma_1,~~\rho_2 = [\rho_1, \rho_2] \rho_1.
$$
then
$\overline{\sigma}_2 = \overline{\sigma}_1$ and $\overline{\rho}_2 = \overline{\rho}_1$ in $G$.
So $G = \langle \overline{\sigma}_1,  \overline{\rho}_1\rangle$ is 2-generated and 2-step nilpotent
hence its commutator subgroup $\Gamma_2(G)$ is cyclic. In $\Gamma_2(G)$ the following relation is true
$[\overline{\sigma}_1, \overline{\rho}_1]^2 = [\overline{\sigma}_1, \overline{\rho}_1^2] =1.$
Hence the commutator subgroup $\Gamma_2(G)$ is generated by $[\overline{\sigma}_1, \overline{\rho}_1]$ and
has order $\leq 2.$

To see that $[\overline{\sigma}_1, \overline{\rho}_1] \neq 1$ we recall that  the unitriangular
group $UT_3(\mathbb{Z})$ over $\mathbb{Z}$ is generated by two transvections, $t_{12}(1)$ and
$t_{23}(1)$ (where  $t_{ij}(1) = e + e_{ij}$), and it is a free 2-step nilpotent group.

Now, there is a
homomorphism $\varphi$
of $G$ onto the unitriangular group
$UT_3(\mathbb{Z}_2),$ where $\mathbb{Z}_2 = \mathbb{Z}/2\mathbb{Z} = \{\overline{0}, \overline{1} \}$,
 by the
rule $\varphi(\overline{\sigma}_1) = t_{12}(\overline{1}),$
$\varphi(\overline{\rho}_1) = t_{23}(\overline{1}),$ It is easy
to see that
$$
\varphi([\overline{\sigma}_1, \overline{\rho}_1])=[t_{12}(\overline{1}),
t_{23}(\overline{1})] = t_{13}(\overline{1}) \neq e.
$$
Hence, $[\overline{\sigma}_1, \overline{\rho}_1] \neq 1$ in $G$ and so
$\Gamma_2(VB_3)/\Gamma_3(VB_3) = \Gamma_2(G)$ is isomorphic to $\Z_2$.

Now, in order to prove the  statement for $n\ge4$, one can   easily adapt an argument
proposed in~\cite{BGG}.
Denote $\Gamma_i=\Gamma_i(VB_n)$ and consider the following short exact sequence:
\begin{equation*}
1 \to \Gamma_2 / \Gamma_3 \to
  \Gamma_1 / \Gamma_3 \stackrel{p}{\to} \Gamma_1/\Gamma_2\to 1,
\end{equation*}
Since any generator $\sigma_i$ in $\Gamma_1/\Gamma_3$ projects to the same
element of $\Gamma_1/\Gamma_2$, for each $1\leq i\leq n-1$, there
exists $t_i\in \Gamma_2/\Gamma_3$ (with $t_1=1$) such that
$\sigma_i=t_i\sigma_1$. Projecting the braid
relation~(\ref{eq2})  into $\Gamma_1/\Gamma_3$, we see that
$t_i \sigma_1 t_{i+1}\sigma_1 t_i \sigma_1= t_{i+1} \sigma_1
t_i\sigma_1 t_{i+1} \sigma_1$. But the $t_i$ are central in $\Gamma_1 /
\Gamma_3$, so $t_i=t_{i+1}$, and since $t_1=1$, we obtain
$\sigma_1=\ldots =\sigma_{n-1}$. In the same way one obtains that
$\rho_1=\ldots =\rho_{n-1}$ in $\Gamma_1/\Gamma_3$. From
relation~(\ref{eq20}) one deduces that
$\rho_1$ and $\sigma_{1}$ commute in $\Gamma_1/\Gamma_3$. Therefore, the surjective
homomorphism $p$ is
in fact an isomorphism.
\item The group $VB_n$ contains $B_n$ (see Proposition~\ref{embeds}) which is not residually
nilpotent for $n\ge3$.
\end{enumerate}
\end{proof}

The group  $\Gamma_2(B_n)$ is perfect (i.e. $\Gamma_2(B_n)=[\Gamma_2(B_n),\Gamma_2(B_n)]$)
for $n\ge 5$~\cite{GoL}. 
An analogous result holds for the group $\Gamma_2(VB_n)$.

\begin{prop}
The group $\Gamma_2(VB_n)$ is perfect for $n\ge 5$.
\end{prop}
\begin{proof}
Let $A_n$ be the alternating group.
The groups  $\Gamma_2(B_n)$ and $\Gamma_2(S_n)=A_n$ are perfect for
$n \geq 5$ (moreover, $A_n$  is simple  for $n \geq 5$). 
Consider a commutator $[u, v] \in \Gamma_2(VB_n)$. In order to prove the claim
we need to show  that $[u, v] \in \Gamma_3(VB_n)$. Since 
$VB_n = \langle B_n, S_n \rangle$ we can use the following commutator identities
$$
[ab, c] = [a, c]^b [b, c],~~~[a, b c] = [a, c] [a, b]^c,~~~[a, b] = [b, a]^{-1},~~~
[a^{-1}, b] = [b, a]^{a^{-1}},
$$
and we can represent the commutator $[u, v]$ as a product of commutators
$$
[\sigma_i, \sigma_j]^{\alpha},~~~[\rho_i, \rho_j]^{\beta},~~~[\sigma_i, \rho_j]^{\gamma},~~~
1 \leq i, j \leq n-1,~~~\alpha, \beta, \gamma \in VB_n.
$$
Since $\Gamma_2(B_n)$ and $\Gamma_2(S_n)$ are perfect for  $n \geq 5$, then $[\sigma_i, \sigma_j] 
\in [\Gamma_2(B_n),\Gamma_2(B_n)]$,
$[\rho_i, \rho_j] \in [\Gamma_2(S_n),\Gamma_2(S_n)]$, and so $[\sigma_i, \sigma_j]^{\alpha}$ and
$[\rho_i, \rho_j]^{\beta}$ belong to  $[\Gamma_2(VB_n),\Gamma_2(VB_n)]$. Therefore we need only to prove that 
commutators of type $[\sigma_i, \rho_j]$ belong to  $[\Gamma_2(VB_n),\Gamma_2(VB_n)]$.

Consider $[\sigma_i, \rho_j]$. If $|i - j| > 1,$ then $[\sigma_i, \rho_j] = 1$.
Let $|i - j| \leq 1.$ Then there are  a pair $k, l$, $1\leq k, l \leq n-1,$ $|k-l| > 1$
and two elements  $c_{i,k} \in \Gamma_2(B_n)$ and  $d_{j,l} \in \Gamma_2(S_n)$ such that $\sigma_i = c_{i,k} \sigma_k,$
$\rho_j = d_{j,l} \rho_l.$ Hence,
$$
[\sigma_i, \rho_j] = [c_{i,k} \sigma_k, d_{j,l} \rho_l].
$$
Using commutator identities we have
$$
[c_{i,k} \sigma_k, d_{j,l} \rho_l] = [c_{i,k} \sigma_k, \rho_l] [c_{i,k} \sigma_k, d_{j,l}]^{\rho_l}
= [c_{i,k}, \rho_l]^{\sigma_k}  [\sigma_k, \rho_l] \left( [c_{i,k}, d_{j,l}]^{\sigma_k}
[\sigma_k, d_{j,l}] \right)^{\rho_l} =
$$
$$
= [c_{i,k}, \rho_l]^{\sigma_k}  [c_{i,k}, d_{j,l}]^{\sigma_k \rho_l}
[\sigma_k, d_{j,l}]^{\rho_l}.
$$
It is clear that  $[c_{i,k}, d_{j,l}] \in [\Gamma_2(VB_n),\Gamma_2(VB_n)]$ 
and therefore  also $[c_{i,k}, d_{j,l}]^{\sigma_k \rho_l}$ belongs to $[\Gamma_2(VB_n),\Gamma_2(VB_n)]$.
Now, let us consider the commutator
$$
[c_{ik}, \rho_l]^{\sigma_k} = \left( c_{ik}^{-1} c_{ik}^{\rho_l} \right)^{\sigma_k}.
$$
Since $c_{i,k}$ lies in $\Gamma_2(B_n)$ and $\Gamma_2(B_n)$ 
is perfect, then $c_{i,k}$ lies in $[\Gamma_2(B_n),\Gamma_2(B_n)]$ and
hence $[c_{i,k}, \rho_l]^{\sigma_k} \in [\Gamma_2(VB_n),\Gamma_2(VB_n)]$. Analogously one can prove that the commutator
$[\sigma_k, d_{j,l}]^{\rho_l}$ belongs to  $[\Gamma_2(VB_n),\Gamma_2(VB_n)]$. Hence  $\Gamma_2(VB_n)$ is perfect for $n \geq 5$.
\end{proof}

\section{Generators and defining relations of the virtual pure braid group $VP_n$}\label{vp}

In this section we study the virtual pure braid group $VP_n$, introduced in \cite{B}.
Define the map
$$
\nu : VB_n \longrightarrow S_n
$$
of $VB_n$ onto the symmetric group $S_n$ as follows:
$$
\nu(\sigma_i) = \nu(\rho_i) = \rho_i, ~~~ i = 1, 2, \ldots, n-1,
$$
where $S_n$ is generated by $\rho_i$ for $i=1,2,\dots, n-1$.
The kernel $\ker \nu$  is called the
{\it virtual pure braid group on $n$ strands} and it is denoted by $VP_n$.

Define the following elements
$$
\lambda_{i,i+1} = \rho_i \, \sigma_i^{-1},~~~
\lambda_{i+1,i} = \rho_i \, \lambda_{i,i+1} \, \rho_i = \sigma_i^{-1} \, \rho_i,
~~~i=1, 2, \ldots, n-1,
$$
$$
\lambda_{i,j} = \rho_{j-1} \, \rho_{j-2} \ldots \rho_{i+1} \, \lambda_{i,i+1} \, \rho_{i+1}
\ldots \rho_{j-2} \, \rho_{j-1},
$$
$$
\lambda_{j,i} = \rho_{j-1} \, \rho_{j-2} \ldots \rho_{i+1} \, \lambda_{i+1,i} \, \rho_{i+1}
\ldots \rho_{j-2} \, \rho_{j-1}, ~~~1 \leq i < j-1 \leq n-1.
$$

\begin{thm}{\bf(\cite{B})} \label{theorem1}
The group $VP_n$ admits a presentation with the generators $\lambda_{k,\, l},$ $1 \leq k \neq l \leq
n$,
and the defining relations:
\begin{equation}
\lambda_{i,j} \,  \lambda_{k,\, l} = \lambda_{k,\, l}  \, \lambda_{i,j} \, ;
 \label{eq29}
\end{equation}
\begin{equation}
\lambda_{k,i} \,  \lambda_{k,j} \,  \lambda_{i,j} = \lambda_{i,j} \,  \lambda_{k,j} \,  \lambda_{k,i} \, ,
 \label{eq28}
\end{equation}

where  distinct letters stand for distinct indices.
\end{thm}

\begin{rem}\label{YB}
It is worth to remark that the group $VP_n$ has been independently defined and studied
in \cite{BE} in relation to Yang Baxter equations.
More precisely, according to \cite{BE}, the virtual pure braid group on $n$ strands
is called the $n$-th  quasitriangular group $QTr_n$ and it is  generated by
$R_{i,j}$ with $1\le i\not=j \le n$ with defining relations
given by the quantum Yang Baxter equations
\begin{eqnarray*}
R_{i,j} \,  R_{k,\, l} &=& R_{k,\, l}  \, R_{i,j} \, ; \\
R_{k,i} \,  R_{k,j} \,  R_{i,j} &=& R_{i,j} \,  R_{k,j} \,  R_{k,i} \, .
\end{eqnarray*}
\end{rem}

Since the natural section $s: S_n \to VB_n$ is well defined one deduces that $VB_n = VP_n \rtimes S_n$.
Moreover we can characterize the conjugacy action of $S_n$ on $VP_n$.
Let  $VP_n \rtimes S_n$
be the semidirect product defined by the action of $S_n$
on the set
$\{ \lambda_{k,\,l} ~ \vert 1 ~ \leq k \neq l \leq n \}$ by permutation of indices.

\begin{prop}{\bf(\cite{B})} \label{prop1}
The map $\omega: VB_n \to  VP_n \rtimes S_n$  sending any element $v$ of $VB_n$
into $(v ((s\circ \nu)(v))^{-1}, \,\nu(v)) \in VP_n \rtimes S_n$ is an isomorphism.
\end{prop}

Let us define the subgroup
$$
V_i = \langle \lambda_{1,i+1}, \lambda_{2,i+1}, \ldots, \lambda_{i,i+1};
\lambda_{i+1,1}, \lambda_{i+1,2}, \ldots, \lambda_{i+1,i} \rangle, \, i=1, \ldots, n-1,
$$
of $VP_n$.  Let $V_i^*$ be
the normal closure of $V_i$ in $VP_{n}$.
We have a ``forgetting map''
$$
\varphi : VP_n \longrightarrow VP_{n-1}
$$
which takes generators $\lambda_{i,n}$, $i = 1, 2, \ldots , n-1,$
and $\lambda_{n,i}$, $i = 1, 2, \ldots , n-1,$  to the unit and
fixes other generators.
The kernel
of $\varphi$ is the group $V_{n-1}^*$, which turns out to be a free group infinitely
generated.

\begin{thm}{\bf(\cite{B})}  \label{theorem2}
The group $VP_n$, $n \geq 2$, is representable as the semi-direct product
$$
VP_n = V_{n-1}^* \rtimes VP_{n-1} = V_{n-1}^* \rtimes
(V^*_{n-2} \rtimes (\ldots \rtimes
(V_2^* \rtimes V_1^*))\ldots),
$$
where $V_1^*$ is a free group of rank $2$ and $V_i^*$, $i=2,
3, \ldots, n-1,$ are free groups infinitely generated.
\end{thm}

The group $V_{n-1}^*$  is  the normal closure of the set
$\{ \lambda_{1,n}, \lambda_{2,n},$ $\ldots, $ $\lambda_{n-1,n}$,
$\lambda_{n,1}, \lambda_{n,2},$ $ \ldots, \lambda_{n,n-1}\}$.
We refer to \cite{B} for a (infinite) free family of generators.
In the following we set $a^{\pm b}$ for $b^{-1} a^{\pm 1} b$.

\begin{lem}{\bf(\cite{B})}   \label{lemma5}
The following formulae are fulfilled in the group $VP_n$:\\

1)~$ \lambda_{n,\, l}^{\lambda_{i,j}^{\varepsilon}} = \lambda_{n,\, l},
~~~\mbox{max}\{i, j\} < \mbox{max}\{n, l\} \, ,~~\varepsilon = \pm 1;
 $ \\

2)~$\lambda_{i,n}^{\lambda_{i,j}} = \lambda_{n,j}^{\lambda_{i,j}} \lambda_{i,n}
\lambda_{n,j}^{-1},~~\lambda_{i,n}^{\lambda_{i,j}^{-1}} = \lambda_{n,j}^{-1} \lambda_{i,n}
\lambda_{n,j}^{\lambda_{i,j}^{-1}},~i < j < n~\mbox{or}~j < i < n;
 $\\

3)~$ \lambda_{n,i}^{\lambda_{i,j}} = \lambda_{n,j} \lambda_{n,i}
\lambda_{n,j}^{-\lambda_{i,j}},~\lambda_{n,i}^{\lambda_{i,j}^{-1}} = \lambda_{n,j}^{-\lambda_{i,j}^{-1}}
\lambda_{n,i} \lambda_{n,j},~~i < j < n~\mbox{or}~j < i < n;
 $\\

4)~$ \lambda_{j,n}^{\lambda_{i,j}} = \lambda_{i,n} \lambda_{j,n}
\lambda_{n,j} \lambda_{i,n}^{-1} \lambda_{n,j}^{-\lambda_{i,j}},~\lambda_{j,n}^{\lambda_{i,j}^{-1}} = \lambda_{j,n}^{-\lambda_{i,n}^{-1}}
\lambda_{i,j}^{-1} \lambda_{j,n} \lambda_{n,j} \lambda_{i,j},~i < j < n~\mbox{or}~j < i <n,
 $\\

where different letters stand for different indices.
\end{lem}

Lemma~\ref{lemma5} provides the action of   $VP_{n-1}$ on $ V_{n-1}^*
$. Therefore one can deduce that
the semi-direct product given in Theorem~\ref{theorem2} fails to be an
almost-direct product.

Nevertheless we can give a partial result on lower central series of virtual pure braids.
Let us start considering the group $VP_3$. It is generated by the elements
$$
\lambda_{2,1},~~\lambda_{1,2},~~\lambda_{3,1},~~\lambda_{3,2},~~\lambda_{2,3},~~\lambda_{1,3},
$$
and defined by relations
$$
\lambda_{1,2} (\lambda_{1,3} \lambda_{2,3}) = (\lambda_{2,3}
\lambda_{1,3}) \lambda_{1,2},~~~ \lambda_{2,1} (\lambda_{2,3}
\lambda_{1,3}) = (\lambda_{1,3} \lambda_{2,3}) \lambda_{2,1},
$$
$$
\lambda_{3,1} (\lambda_{3,2} \lambda_{1,2}) = (\lambda_{1,2}
\lambda_{3,2}) \lambda_{3,1},~~~ \lambda_{3,2} (\lambda_{3,1}
\lambda_{2,1}) = (\lambda_{2,1} \lambda_{3,1}) \lambda_{3,2},
$$
$$
\lambda_{2,3} (\lambda_{2,1} \lambda_{3,1}) = (\lambda_{3,1}
\lambda_{2,1}) \lambda_{2,3},~~~ \lambda_{1,3} (\lambda_{1,2}
\lambda_{3,2}) = (\lambda_{3,2} \lambda_{1,2}) \lambda_{1,3}.
$$

Define the following order on the generators:
$$
\lambda_{1,2} < \lambda_{2,1} < \lambda_{1,3} < \lambda_{2,3} <
\lambda_{3,1} < \lambda_{3,2}.
$$
With this order we can consider the following  15
commutators as the basic commutators of $VP_3/ \Gamma_3(VP_3)$.
$$
\begin{array}{ll}
& [\lambda_{3,2}, \lambda_{3,1}],  [\lambda_{3,2}, \lambda_{2,3}],
[\lambda_{3,2}, \lambda_{1,3}],
 [\lambda_{3,2}, \lambda_{2,1}], [\lambda_{3,2}, \lambda_{1,2}],\\
 \\
& [ \lambda_{3,1}, \lambda_{2,3}],
[\lambda_{3,1}, \lambda_{1,3}],
 [\lambda_{3,1}, \lambda_{2,1}], [\lambda_{3,1}, \lambda_{1,2}],\\
\\
& [\lambda_{2,3}, \lambda_{1,3}],
 [\lambda_{2,3}, \lambda_{2,1}], [\lambda_{2,3}, \lambda_{1,2}], \\
\\
& [\lambda_{1,3}, \lambda_{2,1}], [\lambda_{1,3}, \lambda_{1,2}],\\
\\
& [\lambda_{2,1}, \lambda_{1,2}]. \\
\end{array}
$$
In $VP_3 / \Gamma_3 (VP_3)$  the defining relations
will have the following form (relations are written in terms of basic commutators):
$$
\begin{array}{ll}
1)~~&  [\lambda_{2,3}, \lambda_{1,3}] [\lambda_{2,3}, \lambda_{1,2}]
[\lambda_{1,3}, \lambda_{1,2}] =1, \\
\\
2)~~& [\lambda_{2,3}, \lambda_{1,3}]^{-1} [\lambda_{1,3}, \lambda_{2,1}]
[\lambda_{2,3}, \lambda_{2,1}] =1,\\
\\
3)~~& [\lambda_{3,2}, \lambda_{3,1}]^{-1} [\lambda_{3,1}, \lambda_{1,2}]
[\lambda_{3,2}, \lambda_{1,2}] =1,\\
\\
4)~~& [\lambda_{3,2}, \lambda_{3,1}] [\lambda_{3,2}, \lambda_{2,1}]
[\lambda_{3,1}, \lambda_{2,1}] =1,\\
\\
5)~~& [\lambda_{3,1}, \lambda_{2,3}]^{-1} [\lambda_{3,1},
\lambda_{2,1}]^{-1} [\lambda_{2,3}, \lambda_{2,1}] =1,\\
\\
6)~~& [\lambda_{3,2}, \lambda_{1,2}] [\lambda_{3,2}, \lambda_{1,3}]
[\lambda_{1,3}, \lambda_{1,2}]^{-1} =1.\\
\end{array}
$$

We see that each from the following commutators
$$
 [\lambda_{2,3}, \lambda_{1,2}],  [\lambda_{1,3}, \lambda_{2,1}],
[\lambda_{3,1}, \lambda_{1,2}],
 [\lambda_{3,2}, \lambda_{2,1}], [\lambda_{3,1}, \lambda_{2,3}], [\lambda_{3,2},
\lambda_{1,3}],
$$
is included only once in the list of relations 1) - 6).
Hence, we can exclude these commutators. Then we get the following result.

\begin{lem}
$\Gamma_2(VP_3) / \Gamma_3(VP_3) \simeq \mathbb{Z}^9$.
\end{lem}

Now, let us outline the general case.
The group $VP_n$ contains $n (n-1)$ generators.
Define the following order on the generators:
$$
\lambda_{1,2} < \lambda_{2,1} < \lambda_{1,3} < \lambda_{2,3} <
\lambda_{3,1} < \lambda_{3,2} < \ldots <
\lambda_{1,n} < \lambda_{2,n} < \ldots <
$$
$$
< \lambda_{n-1, n} < \lambda_{n,1} <
\lambda_{n,2} < \ldots < \lambda_{n,n-1}.
$$
The group $\Gamma_2(VP_n) / \Gamma_3(VP_n)$ is generated by
$M = n (n-1) (n^2 - n -1)/2$
of the basic
commutators $[\lambda_{i,j}, \lambda_{k,\, l}]$ with $\lambda_{i,j}> \lambda_{k,\, l}$.

There are two types of relations in $VP_n$:
$$
\lambda_{i,j} \lambda_{k,\, l} = \lambda_{k,\, l} \lambda_{i,j},
$$
$$
\lambda_{k,i}  \lambda_{k,j} \lambda_{i,j} =
\lambda_{i,j} \lambda_{k,j} \lambda_{k,i}.
$$
The number of  relations of the first type is equal to
$$
 n (n-1) (n-2) (n-3).
$$

Hence, the number of basic commutators which  are trivial in
$\Gamma_2(VP_n) / \Gamma_3(VP_n)$ is equal to
$$
M_1 =n (n-1) (n-2) (n-3) / 2.
$$
The number of  relations of the second type  is equal to
$$
M_2 = n (n-1) (n-2).
$$
The relation of the second type  modulo $\Gamma_3(VP_n)$ has the form
$$
[\lambda_{i,j}, \lambda_{k,j}] [\lambda_{i,j}, \lambda_{k,i}] [\lambda_{k,j},
\lambda_{k,i}] = 1.
$$

Each commutator in this relation is basic or inverse to basic commutators.
 As in the case $n=3$
one can check that
the commutator $[\lambda_{i,j}, \lambda_{k,i}]$ itself and its inverse
(i.e. $[\lambda_{k,i}, \lambda_{i,j}]$) don't include in other relations. Therefore we can
exclude from basic commutators $[\lambda_{i,j}, \lambda_{k,i}]$ and its inverse. Hence,
the number of non-trivial basic commutators in $\Gamma_2(VP_n) / \Gamma_3(VP_n)$
is equal to $M-M_1-M_2,$ and
we have the following result.

\begin{prop}\label{lcsvpb} The group $\Gamma_2(VP_n) / \Gamma_3(VP_n)$, for $n \geq
3$, is a free abelian group of
rank $n (n-1) (2n - 3) / 2$.
\end{prop}

The group $VP_2$ is free of rank $2$.
Remark also that Proposition~\ref{lcsvpb}
implies that  for $n \geq
3$ the quotient  $D_2(VP_n) / D_3(VP_n)$  is a free abelian group of
rank $n (n-1) (2n - 3) / 2$ since for any group $G$,
$\Gamma_i(G)/\Gamma_{i+1}(G)$ is isomorphic to
$D_i(G)/D_{i+1}(G)$ modulo torsion (see for instance \cite{H}).

\vspace{5pt}

\noindent {\bf Question.} The group $VP_n$ is 
residually torsion-free nilpotent?

\section{The group $H_n$} \label{defk}

In this section we prove another decomposition of $VB_n$ as a semidirect product.
Let
$\mu: VB_n \to S_n$ be the morphism
defined as follows:

$$
\mu(\sigma_i)=1, \; \mu(\rho_i)=\rho_i, \;i=1,2,\dots, n-1\, ,
$$
where $S_n$ is generated by $\rho_i$ for $i=1,2,\dots, n-1$.
Let us denote by $H_n$ the normal closure of $B_n$ in $VB_n$.

It is evident that $\ker \mu$ coincides with $H_n$.
Now, define the following elements:

$x_{i,i+1}=\sigma_i$, $x_{i,j}=\rho_{j-1} \cdots \rho_{i+1} \sigma_i \rho_{i+1} \cdots \rho_{j-1}$
for $1 \le i < j-1 \le n-1$,

$x_{i+1,i}=\rho_i \sigma_i \rho_i$
and $x_{j,i}=\rho_{j-1} \cdots \rho_{i+1} \rho_i \sigma_i \rho_i \rho_{i+1} \cdots \rho_{j-1}$
for $1 \le i < j-1 \le n-1$.

Since the subgroup of $VB_n$ generated by  $\rho_1, \dots, \rho_{n-1}$  is isomorphic to
the symmetric group $S_n$ (Proposition~\ref{embeds}), we can define an action of $S_n=\langle \rho_1, \dots, \rho_{n-1} \rangle$
on the set
$\{x_{i,j} \, , \;  1 \le i \not= j \le n \}$ by permutation of indices, i.~e.,
$x_{i,j}^\rho=x_{\rho(i),\rho(j)}$, $\rho \in S_n.$

\begin{lem}\label{action}
Let $\rho \in S_n$.
The element $\rho \, x_{i,j} \rho^{-1}$ is equivalent to $x_{\rho(i),\rho(j)}$ for  $1\le i \not=j \le n-1$.
\end{lem}
\begin{proof}
It is sufficient to prove the statement only for generators $\rho_k$, for $1\le k\le n-1$.
If $\rho_k \,\not= \rho_i, \rho_{i+1}, \rho_{i-1}$ from relation (6) in Theorem~\ref{presvb}
one obtains that $\rho_k \,\sigma_i \rho_k= \sigma_i=x_{i,i+1}=x_{\rho(i),\rho(i+1)}$.
Otherwise $\rho_{i+1} \,\sigma_i \rho_{i+1}=x_{i,i+2}$,  $\rho_i \,\sigma_i \rho_i=x_{i+1,i}$
or, by relation (7) in Theorem~\ref{presvb}, $\rho_{i-1} \,\sigma_i \rho_{i-1}=\rho_i \,\sigma_{i-1} \rho_i=x_{i-1,i}$.

Let $x_{i,j}=\rho_{j-1} \cdots \rho_{i+1} \sigma_i \rho_{i+1} \cdots \rho_{j-1}$ for $1 \le i < j-1 \le n-1$.
\begin{enumerate}[i)]
\item Let $ k >j$ or $k < i-1$. Then $\rho_{k}$ commute with $x_{i,j}$ and  the claim holds.
\item If $k=j$, by definition,  $\rho_j \,x_{i,j} \rho_j=x_{i,j+1}$. Since $\rho_j(i)=i$ and
$\rho_j(j)=j+1$ the claim holds.
\item If $k=j-1$ one deduces that $\rho_{j-1} \,x_{i,j} \rho_{j-1}=x_{i,j-1}=x_{\rho_{j-1}(i),\rho_{j-1}(j-1)}$.
\item If $i<k<j-1$ it suffices to remark that $\rho_k \rho_{j-1} \cdots \rho_{i+1}=$
 $\rho_{j-1} \cdots \rho_{i+1} \rho_{k+1}$ and that $\rho_{k+1} \rho_{i+1} \cdots \rho_{j-1}=$
$\rho_{i+1} \cdots \rho_{j-1}
\rho_k$. Therefore  $\rho_k \,x_{i,j} \rho_k^{-1}=$
$x_{i,j}=x_{\rho_k(i), \rho_k(j)}$.
\item If $k=i$, we remark first that the equality

  $\rho_i \rho_{j-1} \cdots \rho_{i+1} \sigma_i \rho_{i+1} \cdots \rho_{j-1} \rho_i$
$=\rho_{j-1} \cdots  \rho_i$ $ \rho_{i+1} \sigma_i \rho_{i+1} \rho_i \cdots \rho_{j-1}$

holds in $VB_n$. Applying relations (4) and  (7) of Theorem~\ref{presvb}
 and we obtain the following equality:
\begin{eqnarray*}
 \rho_{i} \rho_{i+1} \sigma_i  \rho_{i+1} \rho_{i} = \sigma_{i+1} \, ,
 \end{eqnarray*}
 and  then
$\rho_i \,x_{i,j} \rho_i=x_{i+1,j}=x_{\rho(i), \rho(i+1)}$.
\item  If $k= i-1$ it suffices to remark that
applying relations (4) and (7) of Theorem~\ref{presvb}
 one obtains the following equalities:
\begin{eqnarray*}
 \rho_{i-1}
\sigma_i  \rho_{i-1} = \rho_{i} \rho_{i} \rho_{i-1}  \sigma_i  \rho_{i-1}=\\
 \rho_{i}\sigma_{i-1}  \rho_i \rho_{i-1} \rho_{i-1}= \rho_{i}  \sigma_{i-1}
\rho_i \, .
 \end{eqnarray*}
\end{enumerate}
The case of $x_{i+1,i}$
 and $x_{j,i}$ for $1 \le i < j-1 \le n-1$ are similar and they are left to the reader.
\end{proof}

\begin{prop}\label{prop2}
The group $H_n$ admits a presentation with the generators $x_{k,\, l},$ $1 \leq k \neq l \leq
n$,
and the defining relations:
\begin{equation} \label{eq40}
x_{i,j} \,  x_{k,\, l} = x_{k,\, l}  \, x_{i,j},
\end{equation}
\begin{equation} \label{eq41}
x_{i,k} \,  x_{k,j} \,  x_{i,k} =  x_{k,j} \,  x_{i,k} \, x_{k,j},
\end{equation}
where  distinct letters stand for distinct indices.
\end{prop}
\begin{proof}
We use the Reidemeister--Schreier method (see, for example, Chapter 2.2 of \cite{KMS}).
The set of elements of $\langle \rho_1, \dots, \rho_{n-1} \rangle$  in normal form:
\begin{eqnarray*}
\Lambda_n = \{ (\rho_{i_1} \rho_{i_1-1} \cdots \rho_{i_1-r_1}) (\rho_{i_2} \rho_{i_2-1} \cdots
\rho_{i_2-r_2})\cdots(\rho_{i_p} \rho_{i_p-1} \cdots \rho_{i_p-r_p})
\,| \\
\, 1\le i_1< i_2 < \cdots < i_p \le n-1, \, 0\le r_j < i_j
\}
\end{eqnarray*}
is a Schreier set of coset representatives of $H_n$ in $VB_n$.
Define the map $^- : VB_n \longrightarrow \Lambda_n$ which takes an element
$w \in VB_n$
into the representative $\overline{w}$ from $\Lambda_n$. The element
$w \overline{w}^{-1}$ belongs to $H_n$.
The group $H_n$ is generated by
$$
s_{\rho, a} = \rho a \cdot (\overline{\rho a})^{-1},
$$
where $\rho$ runs over the set $\Lambda_n$ and $a$ runs over the set of generators of
$VB_n$ (Theorem 2.7 of~\cite{KMS}).

It is easy to establish that $s_{\rho, \rho_i} = e$ for any  $\rho \in \Lambda_n$
and any generator
$\rho_i$ of $VB_n$.

On the other hand,
$$
s_{\rho, \sigma_i} = \rho \, \sigma_i \cdot (\overline{\rho \sigma_i})^{-1}=
\rho \, \sigma_i \cdot (\overline{\rho})^{-1}= \rho \, \sigma_i \rho^{-1}=x_{\rho(i),\rho(i+1)} \,.
$$
It follows  that each generator $s_{\rho, \sigma_i}$
is equal to some
$x_{i,j}$, $1 \leq i \neq j \leq n$.
The inverse statement is also true,
i.~e.,
each element $x_{i,j}$ is equal to some generator $s_{\rho, \sigma_i}$.

To find  defining relations of $H_n$ we define
a rewriting process $\tau $. It allows us to rewrite a word which is written in the generators
of $VB_n$ and presents an element in $H_n$ as a word in the generators of $H_n$.
Let us associate to the reduced word
$$
u = a_1^{\varepsilon_1} \, a_2^{\varepsilon_2} \ldots
a_{\nu}^{\varepsilon_{\nu}},~~~\varepsilon_l = \pm 1,~~~a_l \in
\{\sigma_1, \sigma_2, \ldots, \sigma_{n-1}, \rho_1, \rho_2, \ldots, \rho_{n-1}
\},
$$
the word
$$
\tau(u) = s_{k_1,a_1}^{\varepsilon_1} \,  s_{k_2,a_2}^{\varepsilon_2}
\ldots s_{k_{\nu},a_{\nu}}^{\varepsilon_{\nu}}
$$
in the generators of $H_n$, where $k_j$ is a representative of the $(j-1)$th
initial segment
of the word $u$ if $\varepsilon_j = 1$ and $k_j$ is a representative of the $j$th
initial segment of
the word $u$ if
$\varepsilon_j = -1$.

The group $H_n$ is defined by relations
$$
r_{\mu,\rho} = \tau (\rho  \, r_{\mu} \,  \rho^{-1}),~~~\rho \in
\Lambda_n,
$$
where $r_{\mu}$ is the defining relation of $VB_n$ (Theorem 2.9 of \cite{KMS}).
Denote by
$$
r_1 = \sigma_i  \, \sigma_{i+1}  \, \sigma_i \,  \sigma_{i+1}^{-1} \,  \sigma_i^{-1}
 \, \sigma_{i+1}^{-1}
$$
the first relation of $VB_n$. Then
$$
r_{1,e} = \tau(r_1) = s_{e,\sigma_i} \,  s_{\overline{\sigma_i},\sigma_{i+1}} \,
s_{\overline{\sigma_i \sigma_{i+1}} ,\sigma_{i}} \,
s_{\overline{\sigma_i \sigma_{i+1} \sigma_i \sigma_{i+1}^{-1}} ,\sigma_{i+1}}^{-1} \,
s_{\overline{\sigma_i \sigma_{i+1} \sigma_i \sigma_{i+1}^{-1} \sigma_i^{-1}} ,\sigma_{i}}^{-1} \,
s_{\overline{r_1} ,\sigma_{i+1}}^{-1} =
$$
$$
= s_{e,\sigma_i} \,  s_{e,\sigma_{i+1}} \,
s_{e,\sigma_{i}} \,
s_{e,\sigma_{i+1}}^{-1} \,
s_{e,\sigma_{i}}^{-1} \,
s_{e,\sigma_{i+1}}^{-1} =
x_{i,i+1}  \, x_{i+1,i+2} \,  x_{i,i+1} x_{i+1,i+2}^{-1} \,
 x_{i,i+1}^{-1} \,
 x_{i+1,i+2}^{-1} \, .
$$

Therefore, the following relation
$$
x_{i,i+1}  \, x_{i+1,i+2} \,  x_{i,i+1} = x_{i+1,i+2} \,
 x_{i,i+1} \,
 x_{i+1,i+2}
$$
is fulfilled in $H_n$.
The remaining relations $r_{1,\rho}$, $\rho \in \Lambda_n$, can be obtained from this
relation using conjugation by $\rho$:
$$
r_{1,\rho}  = \rho \, x_{i,i+1}  \, x_{i+1,i+2} \,  x_{i,i+1} x_{i+1,i+2}^{-1} \,
 x_{i,i+1}^{-1} \,
 x_{i+1,i+2}^{-1} \rho^{-1}=
$$
$$
=x_{\rho(i),\rho(i+1)} x_{\rho(i+1),\rho(i+2)} x_{\rho(i),\rho(i+1)}
x_{\rho(i+1),\rho(i+2)}^{-1} x_{\rho(i),\rho(i+1)}^{-1} x_{\rho(i+1),\rho(i+2)}^{-1}
$$
and we obtain  relations (\ref{eq41}).

Let us consider the next relation of $VB_n$:
$$
r_2 = \sigma_i  \, \sigma_j  \, \sigma_i^{-1}  \, \sigma_j^{-1},~~~|i - j| \geq 2.
$$
Applying the rewriting process defined above we obtain   that following equality holds in $H_n$:
$$
r_{2,e} = \tau(r_2) = s_{e,\sigma_i}  \, s_{\overline{\sigma_i},\sigma_{j}} \,
s_{\overline{\sigma_i \sigma_{j} \sigma_i^{-1}} ,\sigma_{i}}^{-1} \,
s_{\overline{r_2} ,\sigma_{j}}^{-1} =
x_{i,i+1}^{-1} \, x_{j,j+1}^{-1} \, x_{i,i+1} \,
x_{j,j+1} \, .
$$
Conjugating this relation by all representatives from $\Lambda_n$ and applying Lemma~\ref{action} as above,
we obtain  relations
 (\ref{eq40}).

Let us prove that only trivial relations follow from all other relations of $VB_n$.
It is evident for
relations (\ref{eq17})--(\ref{eq19}) defining the group $S_n$ because
$s_{\rho,\rho_i} = e$ for all
$\rho \in \Lambda_n$ and $\rho_i$.

Consider the mixed relation (\ref{eq21}) (relation (\ref{eq20})
can be considered  similarly):
$$
r_3 = \sigma_{i+1} \, \rho_i \, \rho_{i+1} \, \sigma_i^{-1} \, \rho_{i+1} \,
\rho_i.
$$
Using the rewriting process, we get
$$
r_{3,e} = \tau(r_3) = s_{e,\sigma_{i+1}} \,
s_{\overline{\sigma_{i+1} \rho_i \rho_{i+1} \sigma_i^{-1}}
,\sigma_{i}}^{-1} =
$$
$$
= x_{i+1,i+2} (\rho_i  \, \rho_{i+1} \, x_{i,i+1}^{-1} \, \rho_{i+1} \, \rho_i)  = e.
$$.
\end{proof}

The following Corollary is a  straigthforward consequence of Lemma~\ref{action}
and of the fact that the natural section
$S_n \to VB_n$ is well defined.

\begin{cor}\label{hn}
The group  $VB_n$ is isomorphic to  $H_n \rtimes S_n$ where $S_n$ acts by permutation of indices.
\end{cor}

From Proposition~\ref{prop2} we derive that $H_n$, which is the normal closure of $B_n$, is an Artin-Tits group, but  not of spherical type.
It is also easy to verify that the lower central series of $H_n$ is  similar to the one of $B_n$.

\begin{prop}\label{prop:lcsHn}
The following properties hold:
\begin{enumerate}[a)]
\item The group $H_2$ is isomorphic to $\Z*\Z$  which is
residually nilpotent.
\item The quotient $\Gamma_1(H_3)/\Gamma_2(H_3)$ is isomorphic to $\Z\oplus\Z$ and
if $n\ge4$ then  $\Gamma_1(H_n)/\Gamma_2(H_n)$ is isomorphic to $\Z$.
\item If $n\ge3$ then the group  $H_n$ is not residually nilpotent and $\Gamma_2(H_n)=\Gamma_3(H_n)$.
\item If $n\ge5$ the group  $\Gamma_2(H_n)$ is perfect.
\end{enumerate}
\end{prop}
\begin{proof}
We prove the point b).  If $n=3$, from the six
defining relations
one deduces  that
$x_{1,2}=x_{2,3}=x_{3,1}$ and $x_{1,3}=x_{3,2}=x_{2,1}$ in
$\Gamma_1(H_3)/\Gamma_2(H_3)$ which turns to be isomorphic to $\Z\oplus\Z$.
When $n\ge 4$, given two elements $x_{i,j}$ and $x_{k,\, l}$ we have the
following cases:
\begin{enumerate}[i)]
\item  If $j=k$ and $i \not=l$ from the relation $x_{i,j} x_{j,\, l} x_{i,j}=
  x_{j,\, l} x_{i,j} x_{j,\, l}$ one deduces that  $x_{i,j}=x_{j,\, l} $ in
  $\Gamma_1(H_n)/\Gamma_2(H_n)$.
\item  If $j\not=k$ and $i =l$ we conclude as above that  $x_{k,i}=x_{i,j}$  in
  $\Gamma_1(H_n)/\Gamma_2(H_n)$.
\item If $j\not=l$ and $i =k$ there
  exists $1\le m  \le n$ distinct from $j,\, l,k$ such that
$x_{k,j}= x_{j,m}=x_{m,k}=x_{k,\, l}$ in
  $\Gamma_1(H_n)/\Gamma_2(H_n)$.
\item If $i\not=k$
and $j =l$ we proceed as in previous case and we obtain
  that $x_{i,j}=x_{k,j}$
   in  $\Gamma_1(H_n)/\Gamma_2(H_n)$.
\item If $i,j,k,l$ are distinct, using the element $x_{j,k}$ it is clear that
  $x_{i,j}=x_{k,\, l}$ in   $\Gamma_1(H_n)/\Gamma_2(H_n)$.
\item Finally if $i=l$ and $j=k$, we choose $1\le m,\, p \le n$ distinct from $i$
  and $j$ and we obtain the following sequence of identities
$$ x_{i,j}=x_{j,m}=x_{m,p}=x_{p,j}=x_{j,i}$$
holds in
  $\Gamma_1(H_n)/\Gamma_2(H_n)$.
\end{enumerate}
Therefore all $x_{i,j}$ are identified in  $\Gamma_1(H_n)/\Gamma_2(H_n)$.

To prove c) and d) we recall that  $H_n$ is the normal closure of $B_n$
and therefore is not residually nilpotent for $n\ge 3$.
Moreover, we recall that from the
Artin braid relations, it follows that $\Gamma_2(B_n)$ is the normal
closure in $B_n$ of the element $\sigma_1 \sigma_2^{-1}$ (see for instance \cite{BGG}), and thus
$\Gamma_2(H_n)$ coincides with the normal closure in
$VB_n$ of the element $\sigma_1 \sigma_2^{-1}$. Since $\sigma_1 \sigma_2^{-1}=[[\sigma_1,\sigma_2],\sigma_1]^{\sigma_1}$
in $B_n$, then $\sigma_1 \sigma_2^{-1}=[[\sigma_1,\sigma_2],\sigma_1]^{\sigma_1}$ in $VB_n$ and therefore
$\Gamma_2(H_n)=\Gamma_3(H_n)$. In the same way, since
$\Gamma_2(B_n)$ is perfect for all $n \ge 5$~\cite{GL}, so is $\Gamma_2(H_n)$.
\end{proof}

\begin{rem}
The group  $H_3$ decomposes as  a free product $G_1 * G_2$,
where $G_1 = \langle x_{1,2}, x_{2,3}, x_{3,1} \rangle$,
$G_2 = \langle x_{1,3}, x_{3,2}, x_{2,1} \rangle$. The group $G_i,$ $i=1,2$ is isomorphic to the
$2$nd affine Artin-Tits group of type $\mathcal{A}$, also called \emph{circular
braid group} on $3$ strands (see \cite{AJ}).
\end{rem}

\vspace*{5pt}

The decomposition of $VB_n$ into  semidirect product provided in Corollary~\ref{hn}  was earlier proposed by Rabenda (\cite{R}).
More precisely let $K_n$  be the (abstract) group with  the following group presentation:
\begin{itemize}
\item Generators: $x_{i,j}$ for $1 \le i\not=j \le n$.
\item Relations: $x_{i,j} \, x_{j,k} \, x_{i,j}= x_{j,k}\,
x_{i,j} \, x_{j,k}$ for $i, j, k$ distinct indices.
\end{itemize}

The symmetric group $S_n$ acts transitively on $K_n$ by permutation of
indices: for any $\sigma$ in $S_n$, $x_{i,j}^\sigma= x_{\sigma(i),
  \sigma(j)}$.
Let $G_n$ be the semi-direct product of $K_n$ and $S_n$ defined by  above action
and let $s_1, \ldots, s_n$ be the generators of $S_n$ considered as generators of $G_n$.

Rabenda defined a map $\phi: G_n \to VB_n$ as follows; $\phi(s_i)=\rho_i$
for $i=1, \dots, n-1$,
$\phi(x_{i,i+1})=\sigma_i$, $\phi(x_{i,j})=\rho_{j-1} \cdots \rho_{i+1} $
$\sigma_i \rho_{i+1} \cdots \rho_{j-1}$ for $1 \le i < j-1 \le n-1$,
$\phi(x_{i+1,i})=\rho_i \sigma_i \rho_i$ and $\phi(x_{j,i})=\rho_{j-1} \cdots$
$ \rho_{i+1} \rho_i $
$\sigma_i \rho_i \rho_{i+1}$ $\cdots \rho_{j-1}$
for $1 \le i < j-1 \le n-1$
 and he outlined a proof of the fact that
the morphism $\phi$ is actually an  isomorphism.


\section{The extended pure braid group $EP_n$}

In this section we determine the relations between the group $H_n$ and the group $VP_n$.

\begin{prop}
The group $H_n$ and $VP_n$ are not isomorphic for $n\ge 3$.
\end{prop}
\begin{proof}
It suffices to remark that the abelianisation of $VP_n$ is isomorphic to
$\Z^{n(n-1)}$ and to compare with the abelianisation of $H_n$ (part b)  of Proposition~\ref{prop:lcsHn}).
\end{proof}

Let $\varphi$ be the map $S_n\to Hom(S_n)$
defined by the action of the symmetric group  on itself by conjugacy.

\begin{prop} \label{prop:pressymm}
The semidirect product
$S_n \rtimes_\varphi S_n$ admits the following group presentation:

\noindent {\bf $\bullet$ Generators:}  $s_i, t_i$, $i = 1, 2, \ldots, n-1$

\noindent {\bf $\bullet$ Relations:}
$$
t_i^2= s_{i}^2=1,~~~i=1,2,\ldots,n-1  \\
$$
$$
s_i s_j=s_js_i,~~~|i-j| \geq 2\\
$$
$$
s_i s_{i+1} s_i= s_{i+1}
s_i s_{i+1},~~~i=1,2,\ldots,n-2 \\
$$
$$
t_i t_j=t_jt_i,~~~|i-j| \geq 2\\
$$
$$
t_i t_{i+1} t_i= t_{i+1}
t_i t_{i+1},~~~i=1,2,\ldots,n-2 \\
$$
$$
t_i s_j t_i=s_j,~~~|i-j| \geq 2 \\
$$
$$
t_i s_i t_i=s_i,~~~i=1,2,\ldots,n-1   \\
$$
$$
t_{i+1} s_i t_{i+1}=
s_{i+1} s_i s_{i+1},~~~i=1,2,\ldots,n-2   \\
$$
$$
t_{i-1} s_i t_{i-1}=
s_{i-1} s_i s_{i-1},~~~i=1,2,\ldots,n-1 $$
\end{prop}

\vspace{5pt}

\noindent In the following we set $M_n$ the semidirect product $S_n \rtimes_\varphi S_n$.
Let $VB_n$ and $M_n$ be provided with the group presentations given respectively in Theorem~\ref{presvb} and in Proposition~\ref{prop:pressymm}
and let $\chi: VB_n \to M_n$ be the morphism defined by $\chi(\sigma_i)=s_i$ and $\chi(\rho_i)=t_i$
for $i = 1, 2, \ldots, n-1$.
We call \emph{extended pure braid group on $n$ strands} the kernel $\ker \chi$ and we denote it by $EP_n$.

Let $\eta_1$ and $\eta_2$ be the maps from $M_n$ to $S_n$ defined respectively
as follows:

$\eta_1(s_i)=1$ and $\eta_1(t_i)=\rho_i$ for $i = 1, 2, \ldots, n-1$ ;

$\eta_2(s_i)=\rho_i$ and $\eta_2(t_i)=\rho_i$ for $i = 1, 2, \ldots, n-1$,

\noindent where $\rho_i$,  for $i = 1, 2, \ldots, n-1$, are  the usual generators of $S_n$.

\begin{prop}\label{relations}
Let $VP_n$ and $H_n$ be provided with the group presentations given respectively in
Theorem~\ref{theorem1} and Proposition~\ref{prop2}.
 \begin{enumerate}[i)]
\item The group $H_n$ coincides with $\ker (\eta_1  \circ \chi)$.
\item The group $VP_n$ coincides with $\ker (\eta_2  \circ \chi)$.
\item The group $H_n \cap VP_n$ coincides with $\ker \chi$.
\end{enumerate}
\end{prop}
\begin{proof}
We recall that $\mu(\sigma_i)=1, \; \mu(\rho_i)=\rho_i$, $\nu(\sigma_i)=\rho_i$
and $\nu(\rho_i)=\rho_i$ for  $i=1,2,\dots, n-1$. Therefore
 $\mu=\eta_1 \circ \chi$ and $\nu=\eta_2 \circ \chi$ and  part i) and ii) follow and then
$\ker \chi \subseteq H_n \cap VP_n$.
Now let $x$ be a non trivial element of $H_n$. Since $H_n$ is the normal closure of $B_n$, the element
$\chi(x)$ belongs to the subgroup generated by $s_1, \dots, s_{n-1}$ which is isomorphic to $S_n$.
Since $\eta_2(s_i)=t_i$ for $i=1,2,\dots, n-1$, it follows that $\eta_2(\chi(x))=1$ if and only if $\chi(x)=1$
and then $\eta_2$ is injective on $\chi(H_n)$. Therefore, if $x$ belongs to $H_n \cap VP_n$ then $x$ belongs to $\ker \chi$.
\end{proof}

There is also another possible definition of $EP_n$ as a generalisation of the classical pure braid group
$P_n$ (actually, it contains properly the normal closure of $P_n$ in $VB_n$).

\begin{prop}
Let $\varepsilon : H_n \to S_n$ be the map defined by
$\varepsilon(x_{i,j})=\varepsilon(x_{j,i})= \rho_i^{\rho_{i+1}
\dots \rho_{j-1}}$ for $1\le i<j \le n$.  The morphism
$\varepsilon$ is well defined,
the group $EP_n$ is isomorphic to $\ker \varepsilon$
and the normal closure of $P_n$ in $VB_n$ is properly included in $EP_n$.
\end{prop}
\begin{proof}
The morphism $\varepsilon$ coincides with the restriction to $H_n$ of the morphism $\nu$ and therefore
$EP_n$ is isomorphic to $\ker \varepsilon$.

Denote by $\langle \langle P_n \rangle \rangle_{VB_n}$ the normal closure of $P_n$ in $VB_n$.
Remark that $\varepsilon(x_{1,2}^2)=1$.
Since $\langle \langle P_n \rangle \rangle_{VB_n}$ is actually the normal closure
of $\sigma_1^2=x_{1,2}^2$,
one deduces that $\langle \langle P_n \rangle \rangle_{VB_n} \subset EP_n$.

On the other hand,
let us consider the following exact sequence:
$$
1\to  \langle \langle P_n \rangle \rangle_{VB_n}  \to VB_n \to HB_n \to 1 \, ,
$$
where  $HB_n$ is the group obtained adding relations $\sigma_i^2=1$ to the group presentation
of $VB_n$. The group $HB_n$ contains elements of  infinite order (for instance, consider
the element $(\sigma_i\rho_i)^2$). Therefore it is not isomorphic to $M_n$
and we deduce that  $\langle \langle P_n \rangle \rangle_{VB_n}$ does not coincide with $EP_n$.
 \end{proof}


\section{Finite type invariants for virtual braids}\label{gpv}

We recall a possible definition of finite type invariants for classical braids, which is the algebraical
version of the usual definition via singular braids (see \cite{pap} or
Section 1.3 and Proposition 2.1 of \cite{GP} in the case of surface braids).

In the following $A$ will denote an abelian group. An invariant of braids
is a set mapping $v: B_n \to A$. Any invariant $v: B_n \to A$ extends by linearity to a morphism of $\Z$-modules
$v: \Z [B_n] \to A$.

Now let $V$ be the two-sided ideal of $\Z [B_n]$ generated by $\{\sigma_i-\sigma_i^{-1}, | \, i=1, \dots, n-1\}$
and let $V^d$ the $d$-th power of $V$. We obtain this way  a filtration
$$
\Z[B_n] \supset V \supset V^2 \supset \dots
$$
that we call \emph{Goussarov-Vassiliev filtration} for $B_n$.

A \emph{finite type (Goussarov-Vassiliev) invariant of degree $d$} is a morphism of \Z-modules $v: \Z [B_n] \to A$
which vanishes on $V^{d+1}$.


 The Goussarov-Vassiliev filtration for $B_n$ corresponds
to the $I$-adic filtration of $P_n$  (i.e., the filtration associated to the augmentation ideal of $P_n$).
 More precisely, we denote by $p$ the canonical projection
of $B_n$
on $S_n$ and we recall that the set section
$s: S_n \to B_n$ sending $\rho_i$ to $\sigma_i$ (for $i=1, \dots, n-1$) determines an
isomorphism of  $\Z$-modules $\Pi: \Z[B_n] \to \Z[P_n] \otimes \Z[S_n]$ defined as
$$
\Pi(\beta)=\beta ((p\circ s)(\beta))^{-1} \otimes p(\beta), \mbox{for}\;  \beta \in B_n
$$

\begin{prop}(Papadima~\cite{pap}). The additive isomorphism $\Pi: \Z[B_n]$ $ \to \Z[P_n] \otimes \Z[S_n]$
sends isomorphically  $V^d$  to $I^d(P_n) \otimes \Z[S_n]$ for all $d\in \N^*$, where $I^d(P_n)$ is the $d$-th power of the
augmentation ideal of $P_n$.
\end{prop}

Gonz\'alez-Meneses and Paris proved a similar proposition for braid groups on closed surfaces (Proposition 2.2 of \cite{GP}).

In the case of virtual braids we can define a new notion of finite type invariant. This notion was introduced in \cite{GPV}
for virtual knots.
Let $J$ be the two-sided ideal of $\Z [VB_n]$ generated by $\{\sigma_i-\rho_i \, | \, i=1, \dots, n-1\}$.

We obtain this way  a filtration
$$
\Z[VB_n] \supset J \supset J^2 \supset \dots
$$
that we call \emph{Goussarov-Polyak-Viro filtration} for $VB_n$.
We will call Goussarov-Polyak-Viro (GPV) invariant of degree $d$  a morphism of \Z-modules $v: \Z [VB_n] \to A$
which vanishes on $J^{d+1}$.

This notion of invariant corresponds to the remark that stating with a virtual knot and replacing finitely many crossings
(positive or negative) we eventuelly get the (virtual) unknot.

On the other hand,
one can also  remark that a GPV invariant restricted to classical braids is a Goussarov-Vassiliev invariant.
In fact, since $(\sigma_i - \rho_i) - (\sigma_i^{-1} - \rho_i) =
(\sigma_i - \sigma_i^{-1})$ one deduces that $J^d \supset V^d$ where $V^d$ is the $d$-th power of $V$, the two-sided ideal of $\Z [VB_n]$ generated by $\{\sigma_i-\sigma_i^{-1}, | \, i=1, \dots, n-1\}$.

The Goussarov-Polyak-Viro filtration for $VB_n$ corresponds
to the $I$-adic filtration of $VP_n$. The map $\omega: VB_n \to VP_n \rtimes S_n$
defined in Proposition~\ref{prop1} determines an isomorphism of $\Z$-algebras
$\Omega: \Z[VB_n] \to \Z[VP_n] \otimes  \Z[S_n]$, where  $\Z[VP_n] \otimes  \Z[S_n]$ carries the natural
structure of $\Z$-algebra induced   by the semi-direct  product  $VP_n \rtimes S_n$.

\begin{prop}
The $\Z$-algebras isomorphism $\Omega: \Z[VB_n] \to \Z[VP_n] \otimes  \Z[S_n]$
sends isomorphically  $J^d$ to $I^d(VP_n)\otimes \Z[S_n]$ for all $d\in \N^*$.
\end{prop}
\begin{proof}
Let us denote $VB_n$ by $B$.
In order to prove the Proposition we need only to verify that:
$$
J^d= BI^dB=BI^d=I^dB\, .
$$
In fact using this equivalences we can repeat word by word the proof of Proposition 2.2 in \cite{GP}
and therefore prove the claim.

Now, since $VP_n$ is a normal subgroup it suffices to prove that $J=BIB$. The inclusion $J \subset BIB$
is obvious. In order to prove the other inclusion we only have to show that $p -1 \in J$ for $p \in VP_n$.
Now let $p_0\in VP_n$, let $m$ denote the number of $\sigma_i^{\pm 1}$ (for $i=1, \ldots,  n-1$) in the word $p_0$ and let $p_k$ the word
obtained  replacing each of the first $k$ letters  $\sigma_i^{\pm 1}$ (for 
$i=1, \ldots,  n-1$) by $\rho_i$ in $p_0$.
The word $p_m$  is the identity  in $VB_n$
and therefore $p_0-1= p_0- \sum_{k=1}^{n-1} (p_k -p_k) - p_n=$
$\sum_{l=0}^{n-1} (p_l - p_{l+1})$ belongs to $J$.
\end{proof}

\end{document}